\documentclass[11pt]{article}
\usepackage[utf8]{inputenc}
\usepackage{amssymb,amsmath}
\usepackage{graphicx,units,mathrsfs,bm}
\usepackage{subfigure}
\usepackage[sort&compress,numbers]{natbib} 
\usepackage{hyperref,xcolor}
\usepackage[paperwidth=8.5in, paperheight=11in, portrait, top=1in, bottom=1.in, left=1.15in, right=1.15in]{geometry}
\usepackage{tikz}
\usetikzlibrary{calc,decorations.markings,positioning}
\allowdisplaybreaks
\newcommand\bt{\raise 2pt \hbox{$\bigtriangledown$}\hskip 1.5pt}

\newtheorem{rmk}{Remark}[section]

\def\mI{{\mathcal{I}}}

\newmuskip\pFqskip
\mathchardef\pFcomma=\mathcode`, 
\newcommand*\pGq[5]{%
  \begingroup
  \begingroup\lccode`~=`,
    \lowercase{\endgroup\def~}{\pFcomma\mkern\pFqskip}%
  \mathcode`,=\string"8000
  {}_{#1}\phi_{#2}\biggl[\genfrac..{0pt}{}{#3}{#4}\bigg|#5\biggr]%
  \endgroup
}

\numberwithin{equation}{section}

\newtheorem{prop}{Proposition}[section]

\newcommand{\RED}[1]{{\color{black} #1}}

\begin{document}
\title{\bf Degenerate Sklyanin algebras, Askey-Wilson polynomials 
and Heun operators}
\author{
Julien Gaboriaud\textsuperscript{$1$}\footnote{
E-mail: julien.gaboriaud@umontreal.ca}~, 
Satoshi Tsujimoto\textsuperscript{$2$}\footnote{
E-mail: tujimoto@i.kyoto-u.ac.jp}~, 
Luc Vinet\textsuperscript{$1$}\footnote{
E-mail: vinet@CRM.UMontreal.CA}~, 
Alexei Zhedanov\textsuperscript{$3$}\footnote{
E-mail: zhedanov@yahoo.com} \\[.5em]
\textsuperscript{$1$}\small~Centre de Recherches Math\'ematiques, 
Universit\'e de Montr\'eal, \\
\small~P.O. Box 6128, Centre-ville Station, Montr\'eal (Qu\'ebec), 
H3C 3J7, Canada.\\[.9em]
\textsuperscript{$2$}\small~Department of Applied Mathematics and 
Physics, Graduate School of Informatics,\\
\small~Kyoto University, Yoshida-Honmachi, Sakyo-ku, Kyoto 606-8501,
Japan\\[.9em]
\textsuperscript{$3$}\small~School of Mathematics, 
Renmin University of China, Beijing, 100872, China
}
\date{\today}
\maketitle

\hrule
\begin{abstract}\noindent 
The $q$-difference equation, the shift and the contiguity relations of the
Askey-Wilson polynomials are cast in the framework of the three and
four-dimensional degenerate Sklyanin algebras $\mathfrak{ska}_3$ and 
$\mathfrak{ska}_4$. It is shown that the $q$-para Racah polynomials
corresponding to a non-conventional truncation of the Askey-Wilson 
polynomials 
form a basis for a finite-dimensional representation of $\mathfrak{ska}_4$. 
The first order Heun operators defined by a degree raising
condition on polynomials 
are shown to form a five-dimensional vector space
that encompasses $\mathfrak{ska}_4$. The most general quadratic expression
in the five basis operators and such that it raises degrees by no more than 
one is identified with the Heun-Askey-Wilson operator.\\[1em]
\end{abstract} 
\hrule


\section{Introduction}

Quite some time ago, it was shown \cite{Gorsky1993,Wiegmann1995} that the
Askey-Wilson difference operator could be realized as a quadratic expression
in the generators of the degenerate Sklyanin algebra (of dimension four). 
A little earlier Kalnins and Miller \cite{Kalnins1989} used symmetry
techniques to derive the orthogonality relation of the Askey-Wilson
polynomials and identified to that end interesting ladder operators. 
Over the years the application of the factorization method
\cite{Bangerezako1999} to these polynomials and the study of their structure
relations \cite{Koornwinder2007} brought attention to related elements. 
More recently advances have been made in the elaboration of the theory of
$q$-Heun operators and the Heun-Askey-Wilson \cite{Baseilhac2019} and
rational \cite{Tsujimoto2019} Heun operators have been identified by
focusing on certain raising properties of their actions on appropriate
spaces of functions. The purpose of this report is to stress the connections
between these topics. 

The paper will develop as follows. In Section \ref{sec:realizations} we
shall introduce three operators involving $q$-shifts that realize a
three-dimensional degenerate Sklyanin algebra $\mathfrak{ska}_3$. 
These operators will be ubiquitous; it will be observed that their linear
combination is diagonal on special Askey-Wilson polynomials in base $q$ and,
following \cite{Gorsky1993}, that the most general quadratic expression 
formed with them
yields the full Askey-Wilson operator in base $q^2$. 
The degenerate Sklyanin algebra $\mathfrak{ska}_4$ obtained by Gorsky 
and Zabrodin will also be introduced. It has $\mathfrak{ska}_3$ as a
subalgebra and will be seen to admit a formal embedding of the 
Askey-Wilson algebra. In Section \ref{sec:contiguity}, it will be seen 
that the contiguity operators introduced by Kalnins and Miller in their
treatment of the Askey-Wilson polynomials all belong to a model of 
the degenerate
Sklyanin algebra $\mathfrak{ska}_4$. It will also be seen that the 
$q$-para Racah polynomials \cite{Lemay2018} support a finite-dimensional
representation of the four-dimensional degenerate Sklyanin algebra. 
Section \ref{sec:S-Heun} will indicate how the Askey-Wilson bispectral
operators emerge in this context. We shall consider first order
$q$-difference operators and identify the
conditions for such operators to raise by one the degree of polynomials in
the symmetric variable $x=z+z^{-1}$. This will lead to a five dimensional
vector space of operators. A basis will consist of one lowering operator,
two that stabilize polynomials of a given degree and two that are raising
this degree by one. The lowering and stabilizing operators will coincide
with the operators realizing $\mathfrak{ska}_3$ introduced in Section
\ref{sec:realizations}. A combination involving the two raising operators
will give the realization of the fourth generator of $\mathfrak{ska}_4$
beyond those of $\mathfrak{ska}_3$. The relations obeyed by these five
operators will be found in an Appendix. It will further be seen that the
most general quadratic operator in the five basis elements and not raising
the degree by more than one is the Heun-Askey-Wilson operator
\cite{Baseilhac2019a}. This parallels the fact that in bispectral situations
Heun operators could be defined equivalently as raising operators or as
bilinear expressions in the bispectral operators.
As will be indicated in the conclusion this approach is paving the way to
the definition of Sklyanin-like Heun algebras associated to different
degenerations of the Askey-Wilson grid.

\section{Realizations of degenerate Sklyanin algebras and \penalty-10000
Askey-Wilson operators}\label{sec:realizations}

It is well known that quantum algebras can be realized in terms of
$q$-derivatives (see \cite{Floreanini1993, Dobrev1994}). In the case of
$\mathcal{U}_q(\mathfrak{su}(2))$ for example, the commutation relations
\begin{align}\label{uq}
\begin{aligned}{}
    &[\hat{B},\hat{C}] = \frac{\hat{A}^2-\hat{D}^2}{q-q^{-1}}, \qquad
    [\hat{A},\hat{D}]=0, \\
     \hat{A}\hat{B}=q\hat{B}\hat{A}, \quad & \quad
     \hat{B}\hat{D}=q\hat{D}\hat{B}, \qquad 
     \hat{C}\hat{A}=q\hat{A}\hat{C,} \qquad 
     \hat{D}\hat{C}=q\hat{C}\hat{D}
\end{aligned}
\end{align}
are realized \cite{Sklyanin1983}, \cite{Baseilhac2020} by taking
\begin{align}\label{uqreal}
\begin{aligned}{}
    \hat{A}^{(\nu)}=&\,q^{-\nu}T_+, \qquad\quad
    \hat{B}^{(\nu)}=\frac{z}{2(q-q^{-1})}(q^{2\nu}T_--q^{-2\nu}T_+),\\
    \hat{C}&= \frac{2}{(q-q^{-1})z}(T_+-T_-),\qquad\quad
    \hat{D}^{(\nu)}=q^{\nu}T_-,
    \end{aligned}
\end{align}
where in the case of finite dimensional representations $\nu$ is integer or
half-integer (see below) and where $T_+$ and $T_-$ are the $q$-shift operators
that act as follows on functions of $z$:
\begin{align}
 T_+ f(z)=f(qz),\qquad T_- f(z)=f(q^{-1}z).
\end{align}
We shall in the following look at models built with operators of the divided
difference type. 
\subsection{The three-dimensional degenerate Sklyanin algebra
$\mathfrak{ska_3}$ }
Let $p=(a, b, c, \alpha, \beta, \gamma)$ be a set of parameters. The
generalized three-dimensional Sklyanin algebra $\hat{\mathcal{S}}^p$ as defined
in \cite{Iyudu2017} (see also \cite{Chekhov2019}), is given by
three generators $u, v, y$ and the relations:
\begin{equation}\label{eq:sp3}
 uv -avu -\alpha yy=0, \qquad vy-byv-\beta uu, \qquad yu-cuy-\gamma vv=0. 
\end{equation}
Consider the operators
\begin{align}\label{eq:realizationXYUV}
\begin{aligned}{}
 Y&=\frac{1}{z-z^{-1}}(T_+-T_-),\\
 U&=\frac{1}{z-z^{-1}}(zT_+-z^{-1}T_-),\\
 V&=\frac{1}{z-z^{-1}}(zT_--z^{-1}T_+).
\end{aligned}
\end{align}
It is readily checked that they satisfy the following relations:
\begin{align}\label{relA}
\begin{aligned}
 &VY-qYV=0,\\
 &YU-qUY=0,\\
 &[U,V]=(q-q^{-1})Y^{2}.
\end{aligned}
\end{align}
Under the correspondence $\it{lower case} \rightarrow \it{upper case}$, it is
seen that $Y$, $U$, $V$ realize a special case of $\hat{\mathcal{S}}^p$ 
with 
\begin{equation}\label{eq:paramsak3}
    a=1,\qquad b=c=q,\qquad \alpha= (q-q^{-1}), \qquad \beta= \gamma=0.
\end{equation}
We shall henceforth denote this algebra by $\mathfrak{ska}_3$. As shown in
\cite{Iyudu2017}, it corresponds to one of the situations ($(a, b, c) 
\neq (0, 0, 0) \; \text{and} \; \beta = \gamma = b - c =0)$) for which the
generalized Sklyanin algebra $\hat{\mathcal{S}}^p$ has a polynomial growth
Hilbert series (PHS) and Koszul properties.
The algebra $\mathfrak{ska}_3$ thus defined possesses a quadratic Casimir
elements $\Omega^{(2)}$:
\begin{align}
 \Omega^{(2)}=uv+q^{-1}y^2
\end{align}
that takes the value $1$ in the realization \eqref{eq:realizationXYUV} which
implies that $UV$ is related to $Y^2$.

\subsection{The Askey-Wilson polynomials and algebra}

Let us recall that the Askey-Wilson polynomials $p_n(x; a, b, c, d | q)$
defined by
  \begin{align}
  \begin{aligned} \label{AWPolynomials}
  \frac{a^n p_n(x;a,b,c,d|q)}{(ab, ac, ad; q)_n}
  =\pGq{4}{3}{q^{-n}, abcdq^{n-1}, az, az^{-1}} {ab, ac, ad}{q;q}
  \end{aligned}
  \end{align}
with $x = z+z^{-1}$ are eigenfunctions of the operator 
$\mathcal{L}_{q}^{(a, b, c, d)}$ \cite{Koekoek2010}
\begin{equation}
  \mathcal{L}_{q}^{(a, b, c, d)} p_n(x; a, b, c, d | q)
  =\lambda _n p_n(x; a, b, c, d | q) 
\end{equation}
with eigenvalues
\begin{equation}\label{eigen}
    \lambda_n = q^{-n}(1-q^n)(1-abcdq^{n-1}).
\end{equation}
We use standard notation for the basic hypergeometric functions and $q$-shifted
factorials \cite{Koekoek2010}. In base $q^r$, the Askey-Wilson operator reads
\begin{equation}
    \mathcal{L}_{q^r}^{(a, b, c, d)} 
    = A^{(r)}(z)T_+^r-[A^{(r)}(z) + A^{(r)}(z^{-1})]\mI+A^{(r)}(z^{-1})T_-^{r}
\end{equation}
with
\begin{equation}\label{A}
    A^{(r)}(z) = \frac{(1 -az)(1-bz)(1-cz)(1-dz)}{(1-z^2)(1-q^rz^2)}
\end{equation}
and where $\mI$ is the identity operator.

The Askey-Wilson algebra $AW(3)$ \cite{Zhedanov1991} that encodes the
bispectrality of the polynomials $p_n$ is realized by taking the generators
$K_0 = \mathcal{L}_q^{(a, b, c, d)} + (1+q^{-1}abcd)$ and $K_1=x$ to find 
that the defining relations of $AW(3)$
 \begin{align}\label{ZhedanovAW}
 [K_0,K_1]_q=K_2,\qquad 
\begin{aligned}{}
 [K_1,K_2]_q&=\mu K_1+\nu_0 K_0+\rho_0,\\
 [K_2,K_0]_q&=\mu K_0+\nu_1 K_1+\rho_1,
\end{aligned}
\end{align}
where {$[A,B]_q=q^{1/2}AB-q^{-1/2}BA$}, are verified with the parameters 
$\mu$, $\nu$ and $\rho$ related to those, $a$, $b$, $c$, $d$ of the
polynomials $p_n$ (see for instance \cite{Koornwinder2006}).\\[1em]
Consider now the following general linear combination of $Y$, $U$ and $V$:
\begin{equation}
    \mathcal{M}^{(\alpha,\beta,\gamma)} = \alpha Y + \beta U + \gamma V.
\end{equation}
Using \eqref{eq:realizationXYUV}, we see that
\begin{equation}
   \mathcal{M}^{(\alpha,\beta,\gamma)}=F(z)T_++F(z^{-1})T_- 
\end{equation}
where
\begin{equation}\label{F}
    F(z)=\frac{\gamma(1-az)(1-bz)}{1-z^2}
\end{equation}
with
\begin{equation}
    \frac{\alpha}{\gamma}=(a+b), \qquad \frac{\beta}{\gamma}=-ab.
\end{equation}
Since
\begin{equation}
    F(z)+F(z^{-1})=\gamma(1-ab),
\end{equation}
we observe that
\begin{equation}
    \mathcal{M}^{(\alpha,\beta,\gamma)} 
    = \gamma [\mathcal{L}_q^{(a, b, q^{\frac{1}{2}}, -q^{\frac{1}{2}})}
    +(1-ab)].
\end{equation}
It follows that the eigenfunctions of a linear combination of the operators 
$Y, U, V$ such as $\mathcal{M}^{(\alpha,\beta,\gamma)}$ are special
Askey-Wilson polynomials with the property of being ``symmetric" when looked 
at from the dual perspective where variable and degree are exchanged; this 
is because the diagonal term in $\mathcal{M}^{(\alpha,\beta,\gamma)}$ is
constant. 
Correspondingly, following \cite{Koornwinder2006}, by taking
\begin{align}
 K_0 = \frac{1}{\gamma}\mathcal{M}^{(\alpha,\beta,\gamma)}\qquad
 \text{and}\qquad K_1=x
\end{align}
we find that the Askey-Wilson algebra relations \eqref{ZhedanovAW} 
are satisfied with
\begin{align}
\begin{aligned}
&~~\mu =0, \qquad \nu _0=1, \qquad \rho _0=0,\\
\nu _1=-ab&(q-q^{-1})^2,\qquad \rho _1=(1-q^{-1})(a+b)(ab+q).
\end{aligned}\end{align}
We shall consider next quadratic expressions in the generators of 
$\mathfrak{ska}_3$.

\subsection{The Askey-Wilson operator and $\mathfrak{ska}_3$}\label{AWandska}

An important observation \cite{Gorsky1993} comes from considering the most
general quadratic expression in the operators $\{Y, U, V\}$ representing
$\mathfrak{ska}_3$. Let us go over this. Define as before another general
linear combination of these operators:
\begin{equation}
    \mathcal{M}^{(\delta,\epsilon,\zeta)} 
    = \delta Y + \epsilon U + \zeta V = G(z)T_++G(z^{-1})T_-
\end{equation}
where
\begin{equation}
    G(z)=\frac{\zeta(1-q^{-1}cz)(1-q^{-1}dz)}{1-z^2}
\end{equation}
with
\begin{equation}\label{G}
 \frac{\delta}{\zeta}=q^{-1}(c+d), \qquad \frac{\epsilon}{\zeta}=-q^{-2}cd.
\end{equation}
The product 
$\mathcal{M}^{(\alpha,\beta,\gamma)}\mathcal{M}^{(\delta,\epsilon,\zeta)}$ 
will take the form:
\begin{equation}
\begin{split}
 &\mathcal{M}^{(\alpha,\beta,\gamma)}\mathcal{M}^{(\delta,\epsilon,\zeta)}\\
 &\qquad =F(z)G(qz)T_+^2 +[F(z)G(q^{-1}z^{-1})+F(z^{-1})G(q^{-1}z)]\mI
  +F(z^{-1})G(qz^{-1})T_-^2, 
 \end{split}
\end{equation}
A straightforward computation shows that for the specific functions $F(z)$
and $G(z)$ given in \eqref{F} and \eqref{G}, the following identity holds:
\begin{equation}
    F(z)G(q^{-1}z^{-1})+F(z^{-1})G(q^{-1}z)
    =-F(z)G(qz)-F(z^{-1})G(qz^{-1})+\Gamma 
\end{equation}
with $\Gamma$ a constant given by
\begin{equation}
 \Gamma=\gamma \zeta (abcdq^{-2}-ab-cdq^{-2}+1).
\end{equation}
Recalling the expression of $A^{(2)}(z)$ in \eqref{A}, we see that
\begin{equation}
    F(z)G(qz)=\gamma \zeta A^{(2)}(z)
\end{equation}
and hence we write
\begin{equation}\label{factAW}
   \mathcal{M}^{(\alpha,\beta,\gamma)}\mathcal{M}^{(\delta,\epsilon,\zeta)} 
   = \gamma \zeta \: [\mathcal{L}_{q^2}^{(a, b, c, d)}  
   + (abcdq^{-2}-ab-cdq^{-2}+1)\mI].
\end{equation}
We have thus obtained a factorization of the Askey-Wilson operator 
$\mathcal{L}_{q^2}^{(a, b, c, d)}$ as a product of two linear combinations 
of the generators in the representation \eqref{eq:realizationXYUV} of the
special generalized Sklyanin algebra $\mathfrak{ska}_3$.

We also note that 
\begin{equation}
    \mathcal{M}^{(\alpha,\beta,\gamma)}\mathcal{M}^{(\delta,\epsilon,\zeta)} 
    = (\alpha Y + \beta U + \gamma V)(\delta Y + \epsilon U + \zeta V)
\end{equation}
provides the most general quadratic expression in the three generators 
$\{Y, U, V\}$. Taking into account the relations \eqref{eq:realizationXYUV}
between the generators and the expression of $UV$ \eqref{UV} (and $VU$) in
terms of $Y^2$ provided by the value of the Casimir, the product
$\mathcal{M}^{(\alpha,\beta,\gamma)}\mathcal{M}^{(\delta,\epsilon,\zeta)}$ 
can be reduced to:
\begin{equation}\label{AWop}
 \begin{split}
  &\mathcal{M}^{(\alpha,\beta,\gamma)}\mathcal{M}^{(\delta,\epsilon,\zeta)}
  \\&= \beta \epsilon U^2 + \gamma \zeta V^2 +
  (\alpha \delta - \beta \zeta q^{-1} - \gamma \epsilon q)Y^2
  +(\alpha \epsilon q + \beta \delta) UY 
  +(\alpha \zeta q^{-1}+\gamma \delta) VY 
  +(\beta \zeta + \gamma \epsilon)\mI.
\end{split}
\end{equation}
We thus recover (with a different parametrization) the result of Gorsky and
Zabrodin \cite{Gorsky1993} according to which the Askey-Wilson $q$-difference
operator is a quadratic expression in the generators of $\mathfrak{ska}_3$.
(As a matter of fact this result is presented in \cite{Gorsky1993} in the
context of the four-dimensional degenerate Sklyanin algebra to which we shall
turn in a moment.)
\begin{rmk}{~}
The idea of obtaining operators of interest, like the Askey-Wilson one, 
as 
quadratic expressions in the generators of fundamental algebras has 
precedents. Of note is the identification of the Askey-Wilson algebra as a 
coideal subalgebra of 
$\mathcal{U}_q(\mathfrak{sl}(2))$\cite{Granovskii1993} and the use of 
the realization \eqref{uq} to obtain the difference operator of the big 
$q$-Jacobi polynomials \cite{Baseilhac2017} as a generator in this 
embedding. 
\end{rmk}
Returning to the factorization formula, since $\gamma$ and $\zeta$ only occur
in the global factor, we may set $\gamma=\zeta=1$. Summing up we thus have:
\begin{equation}\label{fact1}
 \mathcal{M}^{(\alpha,\beta,1)}\mathcal{M}^{(\delta,\epsilon,1)}=
  \mathcal{L}_{q^2}^{(a, b, c, d)}
  + (abcdq^{-2}-ab-cdq^{-2}+1)\mI
\end{equation}
with
\begin{align}
\alpha &=(a+b), &\beta &=-ab,\\
\delta &=q^{-1}(c+d), &\epsilon &=-q^{-2}cd.
\end{align}
With the eigenvalues $\lambda_n$ of $\mathcal{L}_{q^2}^{(a, b, c, d)}$ given 
by \eqref{eigen} with $q$ replaced by $q^2$, it is straightforward to see 
that the Askey-Wilson polynomials with base $q^2$ correspondingly verify
\begin{align}
 \left[\mathcal{M}^{(\alpha,\beta,1)}\mathcal{M}^{(\delta,\epsilon,1)}
 \right]\;p_n(x; a, b, c, d | q^2)  = \rho _n \;p_n(x; a, b, c, d | q^2)
\end{align}
with 
\begin{align}\label{rho}
 \rho _n =  q^{-2n}(1-abq^{2n})(1-cdq^{2n-2}).
\end{align}\par
\begin{rmk}
If we were to consider two linear combinations
$\mathcal{M}^{(\bar{\alpha},\bar{\beta},1)}$ and
$\mathcal{M}^{(\bar{\delta},\bar{\epsilon},1)}$ of $Y$, $U$ and $V$ where the
roles of the pairs of parameters $(a, b)$ and $(c, d)$ are exchanged with
respect to $\mathcal{M}^{(\alpha,\beta,1)}$ and
$\mathcal{M}^{(\delta,\epsilon,1)}$, namely if we were to take
\begin{align}
    \bar{\alpha}&=q^{-1}(a+b)=q^{-1}\alpha,
    &\bar{\beta}&=-q^{-2}ab=q^{-2}\beta,\label{newpar1}\\ 
    \bar{\delta}&=(c+d)=q\delta, 
    &\bar{\epsilon}&=-cd=q^2\epsilon,\label{newpar2}
     \end{align}
we would obtain again a factorization of the Askey-Wilson operator 
$\mathcal{L}_{q^2}^{(a, b, c, d)}$ of similar form
\begin{equation}\label{fact2}
 \mathcal{M}^{(\bar{\delta},\bar{\epsilon},1)}
 \mathcal{M}^{(\bar{\alpha},\bar{\beta},1)}
 =\mathcal{L}_{q^2}^{(a, b, c, d)}  + (abcdq^{-2}-abq^{-2}-cd+1)\mI.
\end{equation}
This is because $ \mathcal{L}_{q^2}^{(a, b, c, d)}$ is invariant under the
permutation of the parameters and the exchanges of the operators $\mathcal{M}$
with the $q$-shifts of the parameters given in \eqref{newpar1} and in
\eqref{newpar2} simply amount to permuting the pairs $(a,b)$ and $(c,d)$ in the
constant term of the rhs of \eqref{fact1}. This is in line with the fact that
\begin{equation}
\mathcal{M}^{(\alpha,\beta,1)}\mathcal{M}^{(\delta,\epsilon,1)} -
  \mathcal{M}^{(\bar{\delta},\bar{\epsilon},1)}
  \mathcal{M}^{(\bar{\alpha},\bar{\beta},1)}  = (-ab+cd)(1-q^{-2})\mI
\end{equation}
as is easily checked using the relations \eqref{relA} as well as \eqref{UV}. 
Conversely, 
\begin{equation}
\mathcal{M}^{(\alpha,\beta,1)}\mathcal{M}^{(\delta,\epsilon,1)} -
  \mathcal{M}^{(\bar{\delta},\bar{\epsilon},1)}
  \mathcal{M}^{(\bar{\alpha},\bar{\beta},1)}  
  = (\beta - q^2\epsilon)(1-q^{-2})\Omega^{(2)}
\end{equation}
with the $\mathcal{M}$s taken as the linear combinations of the generators,
can be seen to package (abstractly) the relations between $y, u$ and $v$ 
\eqref{eq:sp3} with parameters \eqref{eq:paramsak3}.
\end{rmk}

\subsection{The four-dimensional degenerate Sklyanin algebra
$\mathfrak{ska}_4$}

The four-dimensional degenerate Sklyanin algebra $\mathfrak{ska}_4$ was
obtained in \cite{Gorsky1993} as a limit of the elliptic algebra originally
introduced by Sklyanin \cite{Sklyanin1983} (see \cite{Smith1994} 
for a mathematically oriented review). It is presented in terms of four
generators $A$, $B$, $C$, $D$ obeying the following homogeneous quadratic
relations:
\begin{align}\label{dsa}
\begin{aligned}
 DC=qCD,\qquad &CA=qAC,\qquad [A,D]=\frac{(q-q^{-1})^{3}}{4}C^{2},\\
 &~~~[B,C]=\frac{A^{2}-D^{2}}{q-q^{-1}},\\
 AB-qBA=& \,qDB-BD =-\frac{q^{2}-q^{-2}}{4}(DC-CA).
\end{aligned}
\end{align}
This algebra possesses two Casimir elements:
\begin{align}\label{Cas}
 \Omega_0=AD+\frac{(q-q^{-1})^{2}}{4q}C^{2},\qquad
 \Omega_1=\frac{q^{-1}A^{2}+qD^{2}}{(q-q^{-1})^{2}}+BC
  +\frac{q+q^{-1}}{4}C^{2}.
\end{align}
We note that the subalgebra generated by $\{A, C, D\}$ is isomorphic to
$\mathfrak{ska}_3$.

It was observed \cite{Gorsky1993} that the degenerate Sklyanin algebra
contracts to $\mathcal{U}_q(\mathfrak{su}(2))$; indeed, if one sets
$A=\epsilon \hat{A}$, $B=\hat{B}$, $C=\epsilon^2\hat{C}$, 
$D=\epsilon \hat{D}$ and let $\epsilon$ go to zero we see that the relations
\eqref{dsa} reduce to \eqref{uq}.
In keeping with the representation theory \cite{Sklyanin1983} of the
Sklyanin algebra, the finite dimensional representations of its degenerate
version are characterized by an integer or half-integer $\nu$ and are of
dimension $(2\nu+1)$. We know from \cite{Gorsky1993} that these can be
realized by associating $A$, $B$, $C$, $D$ to the following $q$-difference
operators (we shall not distinguish here the abstract algebra element 
from its realization):
\begin{align}\label{realization}
\begin{aligned}
 &A=q^{-\nu}U, \qquad C=\frac{2}{(q-q^{-1})}Y, \qquad D=q^\nu V ,\\
 B=\frac{1}{2(q-q^{-1})(z-z^{-1})}
 &\left[q^{2\nu}(z^{2}T_--z^{-2}T_+)-q^{-2\nu}(z^{2}T_+-z^{-2}T_-)
 -(q+q^{-1})(T_+-T_-)\right],
\end{aligned}
\end{align}
where $U$, $V$, $Y$ are as in \eqref{eq:realizationXYUV}. 
In this realization the Casimir elements $\Omega_0$ and $\Omega_1$ take 
the following values:
\begin{align}
 \Omega_0=1,\qquad \Omega_1=\frac{q^{2\nu+1}+q^{-2\nu-1}}{(q-q^{-1})^{2}}.
\end{align}
In light of \eqref{Cas} and \eqref{realization}, the former relation 
restates the already observed fact that $UV$ is related to $Y^2$, namely that
\begin{equation}\label{UV}
UV=1-q^{-1}Y^2.
\end{equation}
In the realization \eqref{realization}, the contraction from the degenerate
Sklyanin algebra to $\mathcal{U}_q(\mathfrak{su}(2)$ amounts to taking $z$
very large. It is quickly seen that in this limit the divided difference
operators $\{A, B, C,D\}$ given above reduce to the $\{\hat{A}, \hat{B},
\hat{C},\hat{D}\}$ of \eqref{uqreal}. 

Now using the variable $x=z+z^{-1}$, it is readily found that $B$ can be
expressed as 
\begin{equation}\label{B}
 B=\frac{1}{2(q-q^{-1})}\left[q^{-2\nu}(q^{-1}xU-Ux) + q^{2\nu}(qxV-Vx) 
 -(q+q^{-1})Y \right] 
\end{equation}
in terms of the operators \eqref{eq:realizationXYUV} realizing
$\mathfrak{ska}_3$. We shall now indicate how $x$ can be expressed as a
formal power series in terms of $A, B, C, D$ by inverting \eqref{B}. In 
light of the commutation relation $[U,V]=(q-q^{-1})Y^2$ given in
\eqref{eq:realizationXYUV} and the Casimir relation \eqref{UV} we have
\begin{equation}
    UV=1-q^{-1}Y^2 \qquad \text{and} \qquad VU=1-qY^2.
\end{equation}
It follows that $V$ has an inverse $V^{-1}$ given by the formal power series
in $Y$ expressed as follows:
\begin{equation}
    V^{-1}=U(1-qY^2)^{-1}=(1-q^{-1}Y^2)^{-1}U.
\end{equation}
Using the relations $Ux-qxU=-(q-q^{-1})Y$ and
$xV-qVx=q(q-q^{-1})Y$, we arrive at
\begin{align}\label{x}
    x
    &=q^{-2\nu}\left[2B+\left(\frac{q+q^{-1}}{q-q^{-1}}
     +q^{2\nu}-q^{-2\nu}\right)Y\right](1-q^{-4\nu}V^{-1}U)^{-1}V^{-1}.
\end{align}
As indicated before the Askey-Wilson operator 
$\mathcal{L}_{q^2}^{(a, b, c, d)}$ and $x$ generate the Askey-Wilson
algebra. Within the realization in terms of divided difference
operators, we saw that $\mathcal{L}_{q^{2}}^{(a, b, c, d)}$ according
to \eqref{fact1} is obtained as a quadratic expression in the
generators of the subalgebra $\mathfrak{ska}_3$ of $\mathfrak{ska}_4$
and just found as per \eqref{x} that $x$ is in the completion of the
latter algebra. We can therefore assert that the Askey-Wilson algebra
can be formally embedded in this realization of $\mathfrak{ska}_4$.

\section{Contiguity operators of the Askey-Wilson polynomials and the
degenerate Sklyanin algebra}\label{sec:contiguity}

In \cite{Kalnins1989}, Kalnins and Miller presented an elegant derivation of
the weight function of the Askey-Wilson polynomials which is based on
symmetry techniques. We here wish to point out that their approach can
actually be cast in the framework of degenerate Sklyanin algebras. Central to
the treatment in \cite{Kalnins1989} are certain contiguity and ladder
operators that will prove familiar. In order to facilitate comparison with
the original reference we shall adopt essentially the same notation; we shall
however use $q^2$ as the base. 

Kalnins and Miller begin their considerations by observing that the
Askey-Wilson polynomials satisfy the following contiguity relation
\begin{equation}\label{cont1}
 \mu^{(a,b,c,d)}p_n(x; a, b, c, d | q^2)
 = \RED{q^{-n}(1-abq^{2n-2})}
 p_n(x; aq^{-1}, bq^{-1}, cq, dq | q^2)
\end{equation}
if $\mu^{(a,b,c,d)}$ is the following operator:
\begin{equation}
    \mu^{(a,b,c,d)}=\frac{1}{(z-z^{-1})}\big(-z^{-1}(1-aq^{-1}z)
    (1-bq^{-1}z)T_+ + z(1-aq^{-1}z^{-1})(1-bq^{-1}z^{-1})T_-\big).
\end{equation}
It is further observed that
\begin{equation}\label{cont2}
 \mu^{(cq,dq,aq^{-1},bq^{-1})}p_n(x; aq^{-1}, bq^{-1}, cq, dq | q^2)
 =\RED{q^{-n}(1-cdq^{2n})}
 p_n(x; a, b, c, d | q^2).  
\end{equation}
We may proceed from here to derive the weight function by requesting that 
it be such that $\mu^{(cq,dq,aq^{-1},bq^{-1})}$ is the formal adjoint of
$\mu^{(a,b,c,d)}$; this is done in \cite{Kalnins1989}.
Let us focus on the fact that in view of \eqref{cont1} and \eqref{cont2}, 
the Askey-Wilson polynomials are eigenfunctions of
$\mu^{(cq,dq,aq^{-1},bq^{-1})}\mu^{(a,b,c,d)}$, namely,
\begin{align}
  &\left[\mu^{(cq,dq,aq^{-1},bq^{-1})}\mu^{(a,b,c,d)}\right] 
  p_n(x; a, b, c, d | q^2) 
  = \bar{\rho}_n \; p_n(x; a, b, c, d | q^2),\label{KMeq} 
\end{align}
with
\begin{align}
  & \bar{\rho}_n=q^{-2n}(1-cdq^{2n})(1-abq^{2n-2}). \label{KMrho}   
\end{align}
Not surprisingly the factorization of the Askey-Wilson operator that this
eigenvalue equation entails will coincide with the one described in the
preceding section. This is readily established by recognizing that 
\begin{align}
\begin{aligned}
  \mu^{(a,b,c,d)}&=(aq^{-1}+bq^{-1})Y-abq^{-2}U+V \\
  \mu^{(cq,dq,aq^{-1},bq^{-1})}&=(c+d)Y-cdU+V
\end{aligned}\quad
\begin{aligned}
  \hspace{-1em}&=\mathcal{M}^{(\bar{\alpha},\bar{\beta},1)},\\ 
  \hspace{-1em}&=\mathcal{M}^{(\bar{\delta},\bar{\epsilon},1)}.
\end{aligned}
\end{align}
These contiguity operators are thus found to belong to the realization
\eqref{eq:realizationXYUV} of the Sklyanin algebra $\mathfrak{ska}_3$ and we
see that
\begin{equation}
    \mu^{(cq,dq,aq^{-1},bq^{-1})}\mu^{(a,b,c,d)}=
    \mathcal{M}^{(\bar{\delta},\bar{\epsilon},1)}
    \mathcal{M}^{(\bar{\alpha},\bar{\beta},1)},
\end{equation}
with the connection with the Askey-Wilson operator provided by \eqref{fact2};
we note moreover that the eigenvalue $\bar{\rho}_n$ in \eqref{KMrho}
coincides with the expression obtained from \eqref{rho} under the exchange
$(a,b)\leftrightarrow (c, d)$.

Kalnins and Miller consider in addition the lowering operator
$\tau^{(a,b,c,d)}$:
\begin{equation}
    \tau^{(a,b,c,d)}= \frac{1}{z-z^{-1}}(T_+-T_-) = Y ,
\end{equation}
which is nothing else than our operator $Y$ (or $C$).
They proceed to find its adjoint ${\tau^{(a,b,c,d)}} ^*$ which reads:
\begin{align}\label{taustar}
\begin{aligned}{}
&{\tau^{(a,b,c,d)}}^*=\\
  &\frac{q^{-1}}{z-z^{-1}}
  \!\!\left[\frac{(1-az)(1-bz)(1-cz)(1-dz)}{z^2}T_+ 
  -\frac{(1-az^{-1})(1-bz^{-1})(1-cz^{-1})(1-dz^{-1})}{z^{-2}}T_-
  \right]\!\!.
\end{aligned}
\end{align}
These operators act as follows on the Askey-Wilson polynomials:
\begin{align}\label{shift}
\begin{aligned}
 &\tau^{(a,b,c,d)}p_n(x; a, b, c, d | q^2) 
  = \RED{q^n}(1-q^{-2n})(1-abcdq^{2n-2})p_{n-1}(x; aq, bq, cq, dq | q^2), \\
 &{\tau^{(a,b,c,d)}} ^*p_{n-1}(x; aq, bq, cq, dq | q^2)
  =\RED{-q^{-n}} p_n(x; a, b, c, d | q^2).
\end{aligned}
\end{align}
The key point is that ${\tau^{(a,b,c,d)}}^*$ can be expressed as a linear
combination of the generators $A$, $B$, $C$ and $D$ of the degenerate Sklyanin
algebra $\mathfrak{ska}_4$. Let $e_1=(a+b+c+d)$,\penalty-10000 
$e_2=(ab+ac+ad+bc+bd+cd)$, $e_3=abc+abd+acd+bcd$
and $e_4=abcd$ be the elementary symmetric functions in the parameters 
$(a, b, c, d)$, one finds indeed that
\begin{align}\label{taustar2}
  {\tau^{(a,b,c,d)}}^*\!\! = 
  q^{-1}\!\!\left[-e_3(e_4)^{-\frac{1}{4}}A-2(q-q^{-1})(e_4)^{\frac{1}{2}}B
  +\tfrac{(q-q^{-1})}{2}[e_2-(q+q^{-1})(e_4)^{\frac{1}{2}}]C
  +e_1(e_4)^{\frac{1}{4}}D\right]
\end{align}
with 
\begin{equation}\label{trunc}
    q^{-2\nu}=(abcd)^{\frac{1}{2}}.
\end{equation}
We thus observe that the contiguity and raising operators $\mu^{(a,b,c,d)},
\tau^{(a,b,c,d)}$ and  ${\tau^{(a,b,c,d)}}^*$ belong to the realization of 
the degenerate Sklyanin algebra which is hence represented on 
the Askey-Wilson
polynomials.
In general $\nu$ as given by the relation \eqref{trunc} above will not be an
integer or half integer and the corresponding representation extends the
finite-dimensional one discussed in Section \ref{sec:contiguity} to an
infinite-dimensional one. 
\begin{prop}
The operator $\mu^{(a,b,c,d)}$, its adjoint $\mu^{(cq,dq,aq^{-1},bq^{-1})}$,
$\tau^{(a,b,c,d)}$ and ${\tau^{(a,b,c,d)}}^*$ form a basis equivalent to the
set $\{A, B, C, D\}$ as a representation of the degenerate Sklyanin algebra
$\mathfrak{ska}_4$. Their action on the Askey-Wilson polynomials 
$p_n(x; a, b, c, d | q^2)$ is provided by \eqref{cont1}, \eqref{cont2} and
\eqref{shift} respectively. The connection formula \cite{Askey1985},
\cite{Gasper2004} of Askey and Wilson can be used to express these
formulae as combinations of polynomials $p_k(x; a, b, c, d | q^2)$, 
$k=0, 1, ...$, with parameters $a, b, c, d$ fixed, that span the
representation space.
\end{prop}
Imposing that the representation be finite-dimensional amounts to enforcing 
the non-\penalty-10000 conventional truncation condition  
\begin{equation}\label{trunc2}
    (q^2)^{-N+1}=abcd, \qquad N=2\nu+1
\end{equation}
for the Askey-Wilson polynomials with base $q^2$. Quite strikingly this leads
to polynomials called $q$-para Racah polynomials that have been recently
characterized \cite{Lemay2018} and which are in particular orthogonal on a
bilattice composed of two Askey-Wilson grids. We wish to stress this result.
\begin{prop}
The $q$-para Racah polynomials with base $q^2$ support a representation 
of the degenerate Sklyanin algebra of dimension $N=2\nu+1$ with $\nu$ 
integer or half-integer.
\end{prop}
\begin{rmk}
Remarkably the operator ${\tau^{(a,b,c,d)}}^*$ also features centrally in
Koornwinder's study \cite{Koornwinder2007} of the structure
relations of the Askey-Wilson polynomials. These relations amount to raising
and lowering relations where in contradistinction with the shift relations
that we considered above (following Kalnins and Miller), the parameters are
not affected. It is shown in \cite{Koornwinder2007} that such a
structure relation is obtained when ${\tau^{(a,b,c,d)}}^*$ (denoted by $L$ in
\cite{Koornwinder2007} with the factor $q^{-1}$ omitted) acts upon
the Askey-Wilson polynomial $p_n(x; a, b, c, d | q)$ with base $q$. Note that
the shift relation in \eqref{shift} acts on polynomials with base $q^2$. It
is also indicated in \cite{Koornwinder2007} that
$(q-1){\tau^{(a,b,c,d)}}^*=[\mathcal{L}_{q}^{(a, b, c, d)},x]$.
\end{rmk}
\begin{rmk}
It is further recognized in \cite{Koornwinder2007} on the basis of
results of Rains \cite{Rains2006} and Rosengren
\cite{Rosengren2007} that the operator ${\tau^{(a,b,c,d)}}^*$
generates a representation of the degenerate Sklyanin algebra
$\mathfrak{ska}_4$. This is ascertained from the relation
\begin{align}\label{skrel}
 {\tau^{(a,b,ce,de^{-1})}}^* {\tau^{(qa,qb,q^{-1}c,q^{-1}d)}}^*
 ={\tau^{(a,b,c,d)}}^* {\tau^{(qa,qb,q^{-1}ce,q^{-1}de^{-1})}}^*
\end{align}
given in \cite{Koornwinder2007} and easily checked from \eqref{taustar}. As
observed by Koornwinder \cite{Koornwinder2007}, it is the trigonometric
specialization of a formula in \cite{Rains2006} giving the defining relations
of the Sklyanin algebra. We show below how \eqref{skrel} encapsulates 
the relations \eqref{dsa} of $\mathfrak{ska}_4$.
\end{rmk}
Consider the expression \eqref{taustar2} for ${\tau^{(a,b,c,d)}}^* $ as 
a linear combination of the operators $A, B, C, D$.
Substituting \eqref{taustar2} in \eqref{skrel} and multiplying by 
$(e_4)^{\frac{1}{4}}e(1-e)^{-1}(ce-d)^{-1}$, one arrives at
\begin{align}\label{eq:rels02}
\begin{aligned}
 0&=(e_4)^{\frac{3}{4}}(q-q^{-1})(a+b)\left[A^{2}-D^{2}-(q-q^{-1})(BC-CB)\right]\\
  &-2(e_4)^{\frac{1}{2}}(q-q^{-1})ab(AB-qBA)\\
  &-2(e_4)(q-q^{-1})q^{-1}(BD-qDB)\\
  &+(e_4)^{\frac{1}{4}}(a+b)(qab-q^{-1}cd)\left[(AD-DA)
    -\tfrac{1}{4}(q-q^{-1})^{3}C^{2}\right]\\
  &-(e_4)^{\frac{1}{2}}\tfrac{(q-q^{-1})}{2}
    \left( ((a+b)^{2}q^{2}-ab-cd)q^{-1}
    +(e_4)^{\frac{1}{2}}(1+q^{-2}) \right)CD\\
  &+(e_4)^{\frac{1}{2}}\tfrac{(q-q^{-1})}{2}
    \left( ((a+b)^{2}q^{2}-abq^{4}-cd)q^{-2}
    +(e_4)^{\frac{1}{2}}(1+q^{-2})q \right)DC\\
  &-\tfrac{(q-q^{-1})}{2}
    \left( ab(q+q^{-1})(e_4)^{\frac{1}{2}}+
    [e_4(2-q^{-2})+(b^{2}cd+a^{2}cd-a^{2}b^{2}q^{2})] \right)AC\\
  &-\tfrac{(q-q^{-1})}{2}
    \left(-abq(q+q^{-1})(e_4)^{\frac{1}{2}}-q^{-1}
    [e_4(2-q^{2})+(b^{2}cd+a^{2}cd-a^{2}b^{2}q^{2})] \right)CA.
\end{aligned}
\end{align}
We shall illustrate how the defining relations of $\mathfrak{ska}_4$ can be
obtained from \eqref{eq:rels02}. First choose $b=-a$ and $c=0$. The equality
\eqref{eq:rels02} implies
\begin{align}
 CA=qAC.
\end{align}
Substituting this back in \eqref{eq:rels02} and multiplying  
by $(e_4)^{-\frac{1}{4}}$ yields
\begin{align}\label{eq:rels03}
\begin{aligned}
 0&=(e_4)^{\frac{1}{2}}(q-q^{-1})(a+b)\left[A^{2}-D^{2}-(q-q^{-1})(BC-CB)
  \right]\\
  &-2(e_4)^{\frac{1}{4}}(q-q^{-1})ab(AB-qBA)\\
  &+(a+b)(qab-q^{-1}cd)\left[(AD-DA)-\tfrac{1}{4}(q-q^{-1})^{3}C^{2}\right]\\
  &-2(e_4)^{\frac{3}{4}}(q-q^{-1})q^{-1}(BD-qDB)\\
  &-(e_4)^{\frac{1}{4}}\tfrac{(q-q^{-1})}{2}
    \left(((a+b)^{2}q^{2}-ab-cd)q^{-1}
    +(e_4)^{\frac{1}{2}}(1+q^{-2}) \right)CD\\
  &+(e_4)^{\frac{1}{4}}\tfrac{(q-q^{-1})}{2}
    \left(((a+b)^{2}q^{2}-abq^{4}-cd)q^{-2}
    +(e_4)^{\frac{1}{2}}(1+q^{-2})q \right)DC\\
  &+\tfrac{(q-q^{-1})(q^{2}-q^{-2})}{2}((e_4)^{\frac{1}{4}}ab
  -q^{-1}(e_4)^{\frac{3}{4}})CA.
\end{aligned}
\end{align}
Once again, choose $c=0$ for instance. The equality \eqref{eq:rels03} implies
\begin{align}
 AD-DA=\frac{(q-q^{-1})^{3}}{4}C^{2}.
\end{align}
Repeating the same kind of argument, one obtains the other relations
\eqref{dsa} that define $\mathfrak{ska}_4$.

Through the realization that we have 
considered here, we have observed so far that the degenerate Sklyanin algebra
$\mathfrak{ska}_4$ is a
basic structure underneath the theory of Askey-Wilson polynomials. Much 
like a supersymmetric Hamiltonian is the ``square" of supercharges, the
Askey-Wilson operator is quadratic in generators realizing
$\mathfrak{ska}_4$. We also saw that this is intimately connected to the
application of Darboux transformations
or of the factorization method \cite{Kalnins1989},
\cite{Bangerezako1999} to this operator. This approach as we
know is based on the identification of raising operators. It has been
realized recently that raising properties can provide a unifying principle 
in the theory of Heun operators. We next take this angle to revisit the
Heun-Askey-Wilson operator \cite{Baseilhac2019} and sort out the place
occupied by the degenerate Sklyanin algebra in this Heun operator picture.

\section{$S$-Heun operators and the Heun-Askey-Wilson
operator}\label{sec:S-Heun}

The standard Heun operator that defines the ordinary second order differential
equation with four regular singularities \cite{Kristensson2010} has the
property of raising the degree of polynomials by one. It can also be obtained
as a bilinear expression in the bispectral operators of the Jacobi polynomials,
namely, multiplication by the variable and the hypergeometric operator
\cite{Grunbaum2017}. Both viewpoints have been built upon to develop a broad
perspective on operators of Heun type and the algebras they realize. The
tridiagonalization method based on the hypergeometric operator has been
generalized to any bispectral situation and the concept of algebraic Heun
operator \cite{Grunbaum2018} has emerged in this fashion. In a nutshell this
construct amounts to forming the generic bilinear expression in the bispectral
operators. The raising property has been used to arrive at Heun operators
defined on various lattices. In summary, one looks in this case for the most
general second-order operator that raises by one the degree of polynomials on
specified grids. Applied to the Askey-Wilson lattice or polynomials, both
approaches have led equivalently to the Heun-Askey-Wilson operator
\cite{Baseilhac2019}. (The Heun-Racah and Heun-Bannai-Ito operators have
similarly been obtained \cite{Bergeron2020}.) 

Let us mention that the Heun-Askey-Wilson operator has been shown
\cite{Tsujimoto2019} to arise as a degeneration of the one-variable
Ruijsenaars-van Diejen Hamiltonian \cite{VanDiejen1994},
\cite{Ruijsenaars2004}, \cite{Takemura2017}. It has also been found that this
operator can be diagonalized with the help of the algebraic Bethe ansatz
\cite{Baseilhac2019a}. We shall expand this by relating here the
Heun-Askey-Wilson operators to our observations on Sklyanin algebras. To that
end, we shall first focus on determining the most general first order operators
acting on the Askey-Wilson grid that raise the degree of polynomials by one. 
We shall call them special Heun operators or S-Heun operators for short. 
These can be viewed as second order operators without  diagonal terms. 
Indeed if the operator \eqref{S_def} given below is multiplied by $T_+$, 
we readily see that it takes the form of a first order operator 
$A_1(z)T_+^2 + A_2(z)$ on a grid with base $q^2$. Looking for S-Heun operators
is in fact a more basic problem than searching for the generic second order
operator with the raising property as a way of arriving directly at the Heun
operator of Askey-Wilson type. It is hence not surprising that there will be
factorization connotations. This undertaking will reveal that the S-Heun
operators form a five-dimensional space that includes the operators 
$(A, B, C, D)$ realizing $\mathfrak{ska}_4$. We shall further observe that 
the Heun-Askey-Wilson operator has a quadratic expression in terms of these
S-Heun operators.

\subsection{The $S$-Heun operators}\label{subsec:sheun}

Before we apply the raising condition to determine the S-Heun operators that
act through 
$q$-differences on the symmetric variable $x=z+z^{-1}$, for reference, let us
first go over the most simple case of 
first order differential operators that raise by one the degree of polynomials
in the variable $z$. Consider the operator $S$
\begin{align}
 S=F(z)\frac{d}{dz}+G(z)
\end{align}
and demand that $Sp_n(z)=\tilde{p}_{n+1}(z)$ with $p_n$ and $\tilde{p}_n$
polynomials of degree $n$.
It is readily seen that the most general admissible functions $F(z)$ and 
$G(z)$ are
\begin{align}
 F(z)=\alpha_0+\alpha_1z+\alpha_2 z^{2},\qquad G(z)=\beta_0+\beta_1z.
\end{align}
$S$ therefore belongs to a $5$-dimensional vector space with the following
natural basis 
\begin{align}\label{differential}
 L&=\frac{d}{dz},\qquad M_1=1,\qquad M_2=z\frac{d}{dz},\qquad R_1=z,\qquad 
  R_2=z^{2}\frac{d}{dz},
\end{align}
obtained by setting
all coefficients $\alpha_i$ and $\beta_i$ equal to zero except for one. 

These operators can be combined to form the usual finite-dimensional
differential realization of dimension $2j+1$ on monomials $z^{n}$,
$n=0,1,\dots$ of the Lie algebra $\mathfrak{sl}_2$, i.e.:
\begin{align}
 J_0=z\frac{d}{dz}-j=M_2-jM_1,\qquad
 J_+=z^2\frac{d}{dz}-2jz=R_2-2jR_1,\qquad
 J_-=\frac{d}{dz}=L.
\end{align}
This  corresponds to  the $q\to1$ limit of the realization of 
$\mathcal{U}_q(\mathfrak{su}_2)$ given in Section \ref{sec:realizations}. 

Consider now the $q$-difference operator 
\begin{equation}
 S=A_1(z) T_+ + A_2(z) T_- \label{S_def}    
\end{equation}
where
$A_{1,2}(z)$ are functions of $z$. 
Note that these S-Heun operators can be viewed as ``square roots" of the 
general (second order) Heun operators used in \cite{Baseilhac2019}.
Impose again a raising condition on
polynomials $P_n(x(z))$ of degree $n$ in $x(z)=z+z^{-1}$:
\begin{equation}
 S P_n(x(z)) = \Tilde{P}_{n+1}(x(z)) \label{SPP} 
\end{equation}
for all $n=0,1,2,\dots$. 

It is sufficient to check property \eqref{SPP} for the
elementary Askey-Wilson monomials 
\begin{equation}\label{chi} 
\chi_n(z) =z^n + z^{-n},    
\end{equation}
that is to verify that 
\begin{equation}\label{S_chi_n}
S \chi_n(z) = \sum_{k=0}^{n+1}~ a_{nk} \chi_k(z) 
\end{equation}
for some coefficients $a_{nk}$.
Let us look at the action of $S$ on the two Askey-Wilson monomials $\chi_n(x)$
of lowest degrees. For $n=0$, the raising condition reads
\begin{equation}
 A_1(z) + A_2(z) = a_{00} + a_{10} \chi_1(z) \label{S_0}   
\end{equation}
and similarly, for $n=1$ we have
\begin{equation}\label{S_1} 
 A_1(z)(zq+z^{-1}q^{-1})+A_2(z)(zq^{-1}+z^{-1}q) 
 =a_{10} +a_{11} \chi_1(z) + a_{12} \chi_2(z)
\end{equation}
where $a_{00}$, $a_{01}$, $a_{10}$, $a_{11}$, $a_{12}$ are  
arbitrary parameters. 
Evaluating the action of $S$ on the higher degree Askey-Wilson monomials
does not give rise to new parameters: the higher coefficients $a_{nk}$ 
with $n\geq2$ are always expressed in terms of the $a_{0k}$ and $a_{1k}$. 
Hence these $5$ parameters account for all the degrees of freedom 
that the most general S-Heun operator defined on the Askey-Wilson grid
possesses.

Combining \eqref{S_0} and \eqref{S_1}, we find for $A_1(z)$ 
\begin{equation}
 A_1(z) = \frac{\pi_4(z)}{z(1-z^2)(1-q^2)}, \label{A1_expl}   
\end{equation}
where $\pi_4(z)$ is a polynomial of degree four:
\begin{align}\label{pi_4}
\begin{aligned}
\pi_4(z)=(a_{12}q-a_{01}){z}^{4}
 &+(qa_{11}-a_{00}){z}^{3}
 -((1+{q}^{2})a_{01}-qa_{10}){z}^{2}\\
 &+q(a_{11}-qa_{00})z
 +q(a_{12}-qa_{01}). 
\end{aligned}
\end{align}
From the observation that both the lhs and rhs of the system
\eqref{S_0}--\eqref{S_1} are invariant under $z\to z^{-1}$, it follows that
\begin{equation}\label{A2_expl}
A_2(z) = A_1(z^{-1}). 
\end{equation}
This leads to the following proposition.
\begin{prop}
The most general S-Heun operators on the Askey-Wilson grid which are required
by definition to be of the form
\eqref{S_def} and to raise by one the degrees of polynomials in $x=z+z^{-1}$
are specified by the functions $A_{1,2}(z)$ given in 
\eqref{A1_expl}--\eqref{A2_expl}.
\end{prop}
As the operator $S$ depends on $5$ free parameters, it gives rise as in the
differential case to a $5$-dimensional linear space of S-Heun operators. 
A natural basis for this space is formed by
three sets which correspond respectively to lowering, stabilizing and raising
operators:

($i$) Taking $a_{10}=1$ as the only non-zero parameter in \eqref{pi_4}
leads to the operator denoted $L$ 
which decreases the degree of any polynomial in $x(z)$ by
$1$ and changes its parity.

($ii$) Taking either $a_{00}=1$ or $a_{11}=1$ as the only non-zero parameter,
one obtains stabilizing operators, denoted either $M_1$ or $M_2$.
Both preserve the degree as well as the parity of any polynomial in $x(z)$.

($iii$) The choice $a_{01}=1$ and all other parameters equal to $0$ leads to 
the raising operator $R_1$, while the choice $a_{01}=q$, $a_{12}=1$ and 
all other parameters 0 yields the operator $R_2$. 
Both increase by one the degree of any polynomial in $x(z)$ 
and change parity.

For the sake of completeness, we give below the full expressions
of these $5$ operators
\begin{align}\label{eq:LM1M2R1R2}
\begin{aligned}{}
 L&=\frac{1}{q-q^{-1}}\frac{1}{z-z^{-1}}\left(T_+-T_-\right),\\
 M_1&=\frac{1}{q-q^{-1}}\frac{1}{z-z^{-1}}\left((qz+q^{-1}z^{-1})T_-
  -(q^{-1}z+qz^{-1})T_+\right),\\
 M_2&=\frac{1}{q-q^{-1}}\frac{1}{z-z^{-1}}(z+z^{-1})(T_+-T_-),\\
 R_1&=\frac{1}{q-q^{-1}}\frac{1}{z-z^{-1}}(z+z^{-1})
  \left((qz+q^{-1}z^{-1})T_--(q^{-1}z+qz^{-1})T_+\right),\\
 R_2&=\frac{1}{q-q^{-1}}\left(\frac{q^{2}}{z-z^{-1}}(z+z^{-1})(zT_-
  -z^{-1}T_+)-(zT_-+z^{-1}T_+)\right).
\end{aligned}
\end{align}\par
\begin{prop}
The operators $L$, $M_1$, $M_2$, $R_1$, $R_2$ are linearly
independent. They form a basis for the linear space of  
S-Heun operators.
\end{prop}
Note that the $3$ operators $L$, $M_1$, $M_2$ span the $3$-dimensional
subspace of all ``stabilizing"\penalty-10000 S-Heun operators. This means that 
any operator $S=\alpha_0 L + \alpha_1 M_1 + \alpha_2 M_2$ preserves
the degree of any polynomial in $x(z)$, if at least one of $\alpha_1$,
$\alpha_2$ is nonzero.
Comparing \eqref{eq:realizationXYUV} and \eqref{eq:LM1M2R1R2}, it is immediate
to see that
\begin{align}
 Y=(q-q^{-1}) L,\qquad U=\RED{M_1+q M_2},\qquad
 V=\RED{M_1+q^{-1}M_2},
\end{align}
and that the operators $(L, M_1, M_2)$ equivalently realize $\mathfrak{ska}_3$.
We can thus rephrase as follows the observations of Subsection \ref{AWandska}
according to which the Askey-Wilson operator
is given as a quadratic expression in the operators $(Y, U, V)$ representing
$\mathfrak{ska}_3$:
\begin{prop}
The Askey-Wilson operator can be given as the most general quadratic
combination of the S-Heun operators $L$, $M_1$, $M_2$ that stabilize the 
degree of polynomials in $x(z)$.
\end{prop}
We know that the operators $A, C, D$ in the realization \eqref{realization} of
$\mathfrak{ska}_4$  are proportional to $U, Y, V$ respectively. It is not
difficult to see that $B$ in that same realization can be given as the
following combination of $L, R_1, R_2$:
\begin{align}
 B=\frac{(q+q^{-1})[(q^{2\nu}-q^{-2\nu})-(q-q^{-1})]}{2(q-q^{-1})}L
 +\frac{q^{1-2\nu}}{2}R_1
 +\frac{(q^{2\nu-1}-q^{1-2\nu})}{2(q-q^{-1})}R_2.
\end{align}
We thus have:
\begin{prop}
The realization \eqref{realization} of $\mathfrak{ska}_4$ is obtained from
linear combinations of S-Heun operators on the Askey-Wilson grid.
\end{prop}
In addition, the operator $x$ can be constructed as a quadratic polynomial 
in the the elementary S-Heun operators; we have indeed:
\begin{align}
 x=\frac{1}{q^{2}-q^{-2}}\left[(1+q^{-4})(qM_2R_2-R_2M_2)
  +2q^{-3}(qM_1R_2-R_2M_1)\right].
\end{align}
It follows that the Askey-Wilson algebra can be realized by combining
quadratically the five basic S-Heun operators.

\subsection{Heun-Askey-Wilson and S-Heun operators}

We shall now obtain a formula for the Heun-Askey-Wilson operator in terms of
S-Heun operators.

Consider the most general quadratic combination of the operators
$L$, $M_1$, $M_2$, $R_1$, $R_2$ that raises the degree of polynomials in
$x(z)$ by at most one. There should hence be no terms in ${R_1}^2$ and 
${R_2}^2$. Using the relations in the Appendix
\ref{sec:appendix}, one can show that this combination may be written as
follows
\begin{align}\label{eq:Qhawgen}
\begin{aligned}
 Q_{HAW}=&\,\alpha_1{L^{2}}+\alpha_2LM_2+\alpha_3{M_1}^{2}
 +\alpha_4M_1M_2+\alpha_5M_2L+\alpha_6{M_2}^{2}\\
 &+\beta_1M_1R_1+\beta_2R_1M_1+\beta_3R_2M_2
\end{aligned}
\end{align}
where the $\gamma_i$'s and $\delta_i$'s are arbitrary parameters. 
On functions $f(z)$ this operator takes the form:
\begin{align}\label{HAW}
 Q_{HAW}f(z)=[A_1(z)T_+^{2}+A_1(z^{-1})T_-^{2}+A_0(z)\mI]f(z),
\end{align}
with
\begin{align}
 A_1(z)=\frac{Q_6(z)}{z(1-z^{2})(1-q^{2}z^{2})},\qquad
 A_0(z)=-(A_1(z)+A_1(z^{-1}))+p_1(x),
\end{align}
where $Q_6(z)$ is a generic polynomial of degree $6$ in $z$ and
$p_1(x)$ is a generic polynomial of degree $1$ in the variable
$x=z+z^{-1}$. The exact parameters are expressible in terms of 
those of \eqref{eq:Qhawgen}:
\begin{align}
 p_1(x)=\beta_2x+\alpha_3,\qquad\quad
 Q_6(z)=\frac{1}{q^{2}(q-q^{-1})^{2}}\sum_{k=0}^{6}r_kz^{k},\\[.5em]
\begin{aligned}
 r_0=\beta_2 q^4-\beta_3 q^4+\beta_1 q^3+\beta_3 q^2,\qquad\quad
 r_1=\alpha_3 q^4-\alpha_4 q^3+\alpha_6 q^2,\\
 r_2=-\beta_3 q^6+\beta_1 q^5+2\beta_2 q^4+ 
 (\alpha_5+\beta_1)q^3+(\alpha_2+\beta_2-\beta_3)q^2 +\beta_1 q,\\
 r_3=-\alpha_4 q^5+(\alpha_3+\alpha_6)q^4 +\alpha_1 q^3+
 (\alpha_3+\alpha_6)q^2-\alpha_4 q,\\
 r_4=-\beta_3 q^6+\beta_1 q^5+(\alpha_2+\beta_2-\beta_3)q^4 
 +(\alpha_5+\beta_1)q^3 +2 \beta_2 q^2+\beta_1 q,\\
 r_5=\alpha_6 q^4-\alpha_4 q^3+\alpha_3 q^2,\qquad\quad
 r_6=\beta_1 q^3+\beta_2q^2.
\end{aligned}
\end{align}
This operator is recognized as the Askey-Wilson Heun operator which has been
identified and characterized in \cite{Baseilhac2019}.
(See also \cite{Tsujimoto2019} and \cite{Baseilhac2020}.)
It is immediately seen that the Askey-Wilson operator is recovered upon
taking the $\beta_i$'s equal to zero, which is equivalent to removing the
terms involving raising S-Heun operators from $Q_{HAW}$.

This formula giving $Q_{HAW}$ as the most general
quadratic combination in the S-Heun operators on the Askey-Wilson grid
provides a novel characterization of the Heun-Askey-Wilson operator.
As pointed out at the beginning of this section, this operator was identified
in \cite{Baseilhac2019} on the one hand as the most general second order
$q$-shift operator that raises by one the degree of
polynomials on the Askey-Wilson grid and on the other hand, as the
tridiagonalization of the Askey-Wilson operator as per the algebraic Heun
construct. The presentation obtained here with the S-Heun operators as basic
building blocks has the merit of providing, typically, a factorization of $Q_{HAW}$. 
Indeed it is seen that the Heun-Askey-Wilson operator can also be written generically
in the form:
\begin{equation}\label{factQ}
    Q_{HAW}=(\xi_1L+\xi_2M_1+\xi_3M_2)
    (\eta_1L+\eta_2M_1+\eta_3M_2+\eta_4R_1+\eta_5R_2)+\kappa.
\end{equation}
This formula for $Q_{HAW}$ should be compared with equation \eqref{factAW} that
provides the factorization of the Askey-Wilson operator
as the product of two $\mathfrak{ska}_3$ elements. It is hence manifest from
\eqref{factQ} that $Q_{HAW}$ reduces to the Askey-Wilson operator when $g=h=0$.

\section{Conclusion}\label{sec:conclusion}

To conclude, let us first summarize our observations and second offer a brief
outlook.

We have considered realizations of the Sklyanin algebras $\mathfrak{ska}_3$ and
$\mathfrak{ska}_4$ in terms of $q$-difference operators and we determined the
first order operators of that type -- the S-Heun operators -- that are the
basic constituents of the most general degree raising operator in that class.
Within these realizations, our salient observations are:
\begin{itemize}
    \item The Askey-Wilson operator factorizes as the product of two linear
    combinations of elements in $\mathfrak{ska}_3$;
    \item In analogy with the dynamical enlargement of a symmetry algebra with
    the inclusion of ladder operators, the contiguity and shift operators of
    the Askey-Wilson polynomials have been shown to generate a realization of
    the degenerate Sklyanin algebra $\mathfrak{ska}_4$ which formally includes
    a realization of the Askey-Wilson algebra.
    \item The $q$-para Racah polynomials (with base $q^2$) have been identified
    as forming a basis for the finite-dimensional representations of the
    degenerate Sklyanin algebra $\mathfrak{ska}_4$.
    \item The set of S-Heun operators is five-dimensional and has a subset that
    realizes $\mathfrak{ska}_4$.
    \item The operator multiplication by $x$ has a quadratic expression in
    terms of the S-Heun operators.
    \item The Heun-Askey-Wilson operator can also be written as a quadratic
    expression in the S-Heun operators.
\end{itemize}
With respect to these last two points, let us mention the following. We recall
that the algebraic Heun construct gives the Heun-Askey-Wilson operator
$Q_{HAW}$ as a bilinear operator in $x$ and the Askey-Wilson operator. In view
of the first and next to last points, this implies that $Q_{HAW}$ is quartic in
the S-Heun operators, an expression that must be reducible to the quadratic
formula \eqref{eq:Qhawgen} obtained here.

This study raises a number of questions. Let us mention two: ($i$) How does the
examination of the S-Heun operators extend when the raising property is applied
to rational functions as in \cite{Tsujimoto2019}? ($ii$) What other algebraic
structures akin to the degenerate Sklyanin algebras will emerge when the S-Heun
operator approach is adapted to other lattices such as for example the
quadratic one on which the Wilson polynomials are defined?
We plan on addressing these and other related questions in the near future.

\subsection*{Acknowledgments}

The authors benefitted from discussions with Nicolas Cramp\'e and Slava
Spiridonov. JG holds an Alexander-Graham-Bell scholarship from the Natural
Science and Engineering Research Council (NSERC) of Canada.
The work of ST is partially supported by JSPS KAKENHI (Grant 
Numbers 19H01792, 17K18725). The research of LV is funded in part by a
Discovery Grant from NSERC. The work of AZ is supported by the National Science
Foundation of China (Grant No.11771015).

\appendix
\section{Quadratic algebraic relations for $L$, $M_1$, $M_2$, $R_1$, $R_2$}
\label{sec:appendix}
The homogeneous quadratic algebraic relations 
between the five $S$-Heun operators $L$, $M_1$, $M_2$, $R_1$, $R_2$ are
collected below:
\begin{align}
 [M_1,M_2]&=(q+q^{-1})^{2}L^{2}, \label{eqM1M2}\\
 M_1L-(q+q^{-1})LM_1&=LM_2,\\
 LM_1+M_2L&=0,\\
 {M_1}^{2}+{M_2}^{2}+(q+q^{-1})M_2M_1&=1,\\[0.5em]
 LR_1&=1-{M_2}^{2},\\
 R_1L&=1-{M_1}^{2},\\
 LR_2&=-2L^{2}+q^{-1}{M_2}^{2}+M_1M_2+q,\\
 R_2L&=-2L^{2}+q{M_2}^{2}+q^{2}M_2M_1,\\[0.5em]
 R_1M_2+M_1R_1&=0,\\
 M_1R_2+R_2M_2&=2(q+q^{-1})M_2L-(q+q^{-1})^{2}LM_2,\\
 qR_1M_1-M_1R_1&=R_2M_1+(q^{2}+q^{-2})LM_1,\\
 R_1M_1-(q+q^{-1})M_1R_1&=M_2R_1-(q+q^{-1})(q-q^{-1})^{2}M_2L,\\
 M_2R_2-(q+q^{-1})R_2M_2&=R_2M_1+2(q+q^{-1})M_1L
  -(2q^{-2}+1+q^{4})LM_1,\\[0.5em]
 {R_2}^{2}-qR_2R_1+q^{-1}R_1R_2&=-2(q+q^{-1})^{2}M_1M_2
  -(q+q^{-1})^{3}{M_2}^{2}\nonumber\\
  &\quad+2[(q^{2}+q^{-2})
  -(q^{2}-q^{-2})^{2}]L^{2}\label{eqR2R2}
 \end{align}
These relations are checked directly from the expressions of the operators in
\eqref{eq:LM1M2R1R2}. They 
provide the necessary reorderings to
reexpress the the most general quadratic combination of the $5$ operators
as in \eqref{eq:Qhawgen}.

\bibliographystyle{unsrtinur} 
\bibliography{ref_sheun.bib} 

\begin{thebibliography}{10}

\bibitem{Gorsky1993}
A.~S. Gorsky and A.~V. Zabrodin.
\newblock {Degenerations of Sklyanin algebra and Askey-Wilson polynomials}.
\newblock {\em Journal of Physics A: Mathematical and General}, 26(15), 1993.
\newblock \href {http://arxiv.org/abs/hep-th/9303026}
  {\path{arXiv:hep-th/9303026}}.

\bibitem{Wiegmann1995}
P.~B. Wiegmann and A.~V. Zabrodin.
\newblock {Algebraization of difference eigenvalue equations related to
  {$U_q(sl_2)$}}.
\newblock {\em Nuclear Physics, Section B}, 451(3):699--724, 1995.
\newblock \href {http://arxiv.org/abs/cond-mat/9501129}
  {\path{arXiv:cond-mat/9501129}}.

\bibitem{Kalnins1989}
E.~G. Kalnins and W.~J. Miller.
\newblock {Symmetry techniques for {$q$}-series: Askey-Wilson polynomials}.
\newblock {\em Rocky Mountain Journal of Mathematics}, 19(1):223--230, 1989.

\bibitem{Bangerezako1999}
G.~Bangerezako.
\newblock {The factorization method for the Askey–Wilson polynomials}.
\newblock {\em Journal of Computational and Applied Mathematics},
  107(2):219--232, 1999.
\newblock \href {http://arxiv.org/abs/math/9805143}
  {\path{arXiv:math/9805143}}.

\bibitem{Koornwinder2007}
T.~H. Koornwinder.
\newblock {The structure relation for Askey–Wilson polynomials}.
\newblock {\em Journal of Computational and Applied Mathematics}, 207:214--226,
  2007.
\newblock \href {http://arxiv.org/abs/math/0601303}
  {\path{arXiv:math/0601303}}.

\bibitem{Baseilhac2019}
P.~Baseilhac, S.~Tsujimoto, L.~Vinet, and A.~Zhedanov.
\newblock {The Heun-Askey-Wilson Algebra and the Heun Operator of Askey-Wilson
  Type}.
\newblock {\em Annales Henri Poincar{\'{e}}}, 20(9):3091--3112, 2019.
\newblock \href {http://arxiv.org/abs/1811.11407} {\path{arXiv:1811.11407}}.

\bibitem{Tsujimoto2019}
S.~Tsujimoto, L.~Vinet, and A.~Zhedanov.
\newblock {The rational Heun operator and Wilson biorthogonal functions}.
\newblock pages 1--14, 2019.
\newblock \href {http://arxiv.org/abs/1912.11571} {\path{arXiv:1912.11571}}.

\bibitem{Lemay2018}
J.~M. Lemay, L.~Vinet, and A.~Zhedanov.
\newblock {A {$q$}-generalization of the para-Racah polynomials}.
\newblock {\em Journal of Mathematical Analysis and Applications},
  462(1):323--336, 2018.
\newblock \href {http://arxiv.org/abs/1708.03368} {\path{arXiv:1708.03368}}.

\bibitem{Baseilhac2019a}
P.~Baseilhac, L.~Vinet, and A.~Zhedanov.
\newblock {The {$q$}-Heun operator of big {$q$}-Jacobi type and the {$q$}-Heun
  algebra}.
\newblock {\em Ramanujan Journal}, 2019.
\newblock \href {http://arxiv.org/abs/1808.06695} {\path{arXiv:1808.06695}}.

\bibitem{Floreanini1993}
R.~Floreanini and L.~Vinet.
\newblock {$q$}-difference realizations of quantum algebras.
\newblock {\em Physics Letters B}, 315(3-4):299--303, 1993.

\bibitem{Dobrev1994}
V.~K. Dobrev.
\newblock {{$q$}-difference intertwining operators for
  {$U_q(\mathfrak{sl}(n))$}: General setting and the case {$n=3$}}.
\newblock {\em Journal of Physics A: Mathematical and General},
  27(14):4841--4857, 1994.
\newblock \href {http://arxiv.org/abs/hep-th/9405150}
  {\path{arXiv:hep-th/9405150}}.

\bibitem{Sklyanin1983}
E.~K. Sklyanin.
\newblock {Some algebraic structures connected with the Yang-Baxter equation.
  Representations of quantum algebras}.
\newblock {\em Functional Analysis and Its Applications}, 17(4):273--284, 1983.

\bibitem{Baseilhac2020}
P.~Baseilhac, X.~Martin, L.~Vinet, and A.~Zhedanov.
\newblock {Little and big {$q$}-Jacobi polynomials and the Askey–Wilson
  algebra}.
\newblock {\em Ramanujan Journal}, 51(3):629--648, 2020.
\newblock \href {http://arxiv.org/abs/1806.02656} {\path{arXiv:1806.02656}}.

\bibitem{Iyudu2017}
N.~Iyudu and S.~Shkarin.
\newblock {Three dimensional Sklyanin algebras and Gr{\"{o}}bner bases}.
\newblock {\em Journal of Algebra}, 470:379--419, 2017.
\newblock \href {http://arxiv.org/abs/1601.00564} {\path{arXiv:1601.00564}}.

\bibitem{Chekhov2019}
L.~Chekhov, M.~Mazzocco, and V.~Rubtsov.
\newblock {Quantised Painlev\'e monodromy manifolds, Sklyanin and Calabi-Yau
  algebras}.
\newblock pages 1--41, 2019.
\newblock \href {http://arxiv.org/abs/1905.02772} {\path{arXiv:1905.02772}}.

\bibitem{Koekoek2010}
R.~Koekoek, P.~A. Lesky, and R.~F. Swarttouw.
\newblock {\em {Hypergeometric Orthogonal Polynomials and Their
  {$q$}-Analogues}}.
\newblock Springer Monographs in Mathematics. Springer Berlin Heidelberg, 2010.

\bibitem{Zhedanov1991}
A.~S. Zhedanov.
\newblock {“Hidden symmetry” of Askey–Wilson polynomials}.
\newblock {\em Theoretical and Mathematical Physics}, 89(2):1146--1157, 1991.

\bibitem{Koornwinder2006}
T.~H. Koornwinder.
\newblock {The Relationship between Zhedanov's Algebra {$AW(3)$} and the Double
  Affine Hecke Algebra in the Rank One Case}.
\newblock {\em SIGMA. Symmetry, Integrability and Geometry: Methods and
  Applications}, 3:063, dec 2006.
\newblock \href {http://arxiv.org/abs/math/0612730}
  {\path{arXiv:math/0612730}}.

\bibitem{Granovskii1993}
Y.~I. Granovskii and A.~S. Zhedanov.
\newblock {Linear covariance algebra for {$SL_q(2)$}}.
\newblock {\em Journal of Physics A: Mathematical and General}, 26:L357, 1993.

\bibitem{Baseilhac2017}
P.~Baseilhac, L.~Vinet, and A.~Zhedanov.
\newblock {The {$q$}-Onsager algebra and multivariable {$q$}-special
  functions}.
\newblock {\em Journal of Physics A: Mathematical and Theoretical}, 50(39),
  2017.
\newblock \href {http://arxiv.org/abs/1611.09250} {\path{arXiv:1611.09250}}.

\bibitem{Smith1994}
S.~P. Smith.
\newblock {The four-dimensional Sklyanin algebras}.
\newblock {\em K-Theory}, 8(1):65--80, 1994.

\bibitem{Askey1985}
R.~Askey and J.~Wilson.
\newblock {\em {Some basic hypergeometric orthogonal polynomials that
  generalize Jacobi polynomials}}.
\newblock American Mathematical Society, 1985.

\bibitem{Gasper2004}
G.~Gasper and M.~Rahman.
\newblock {\em {Basic Hypergeometric Series}}.
\newblock Cambridge University Press, 2nd edition, 2004.

\bibitem{Rains2006}
E.~M. Rains.
\newblock {{$BC_n$}-symmetric abelian functions}.
\newblock {\em Duke Mathematical Journal}, 135(1):99--180, 2006.
\newblock \href {http://arxiv.org/abs/math/0402113}
  {\path{arXiv:math/0402113}}.

\bibitem{Rosengren2007}
H.~Rosengren.
\newblock {An elementary approach to {$6j$}-symbols (classical, quantum,
  rational, trigonometric, and elliptic)}.
\newblock {\em Ramanujan Journal}, 13(1-3):131--166, 2007.
\newblock \href {http://arxiv.org/abs/math/0312310}
  {\path{arXiv:math/0312310}}.

\bibitem{Kristensson2010}
G.~Kristensson.
\newblock {\em {Second order differential equations: Special functions and
  their classification}}.
\newblock 2010.

\bibitem{Grunbaum2017}
F.~A. {Gr{\"{u}}nbaum}, L.~Vinet, and A.~Zhedanov.
\newblock {Tridiagonalization and the Heun equation}.
\newblock {\em Journal of Mathematical Physics}, 58(3), 2017.
\newblock \href {http://arxiv.org/abs/1602.04840} {\path{arXiv:1602.04840}}.

\bibitem{Grunbaum2018}
F.~A. Gr{\"{u}}nbaum, L.~Vinet, and A.~Zhedanov.
\newblock {Algebraic Heun Operator and Band-Time Limiting}.
\newblock {\em Communications in Mathematical Physics}, 364(3):1041--1068,
  2018.
\newblock \href {http://arxiv.org/abs/1711.07862} {\path{arXiv:1711.07862}}.

\bibitem{Bergeron2020}
G.~Bergeron, N.~Cramp{\'{e}}, S.~Tsujimoto, L.~Vinet, and A.~Zhedanov.
\newblock {The Heun-Racah and Heun-Bannai-Ito algebras}.
\newblock pages 1--18, 2020.
\newblock \href {http://arxiv.org/abs/2003.09558} {\path{arXiv:2003.09558}}.

\bibitem{VanDiejen1994}
J.~F. {Van Diejen}.
\newblock {Integrability of difference Calogero-Moser systems}.
\newblock {\em Journal of Mathematical Physics}, 35(6):2983--3004, 1994.

\bibitem{Ruijsenaars2004}
S.~N.~M. Ruijsenaars.
\newblock {Integrable {$BC_N$} Analytic Difference Operators: Hidden Parameter
  Symmetries and Eigenfunctions}.
\newblock In {\em New Trends in Integrability and Partial Solvability}, pages
  217--261, 2004.

\bibitem{Takemura2017}
K.~Takemura.
\newblock {Degenerations of Ruijsenaars–van Diejen operator and
  {$q$}-Painlev{\'{e}} equations}.
\newblock {\em Journal of Integrable Systems}, 2(1):1--27, 2017.
\newblock \href {http://arxiv.org/abs/1608.07265} {\path{arXiv:1608.07265}}.

\end{thebibliography}

\end{document}